\documentclass[aps,pre,twocolumn]{revtex4-1}
\usepackage{amsmath, amssymb}
\usepackage[pdftex]{graphicx}
\usepackage{color}

\begin{document}

\title{Master Stability Islands for Amplitude Death \\ in Networks of Delay-Coupled Oscillators}

\author{Stanley R. Huddy}
\email{srh@fdu.edu}
\affiliation{Department of Computer Sciences and Engineering, Fairleigh Dickinson University, Teaneck, NJ 07666}

\author{Jie Sun}
\email{sunj@clarkson.edu}
\affiliation{Department of Mathematics, Clarkson University, Potsdam, NY 13699}
\affiliation{Department of Physics, Clarkson University, Potsdam, NY 13699}

\date{\today}

\begin{abstract}
This paper presents a master stability function (MSF) approach for analyzing the stability of amplitude death (AD) in networks of delay-coupled oscillators. Unlike the familiar MSFs for instantaneously coupled networks, which typically have a single input encoding for the effects of the eigenvalues of the network Laplacian matrix, for delay-coupled networks we show that such MSFs generally require two additional inputs: the time delay and the coupling strength. To utilize the MSF for predicting the stability of AD of arbitrary networks for a chosen nonlinear system (node dynamics) and coupling function, we introduce the concept of master stability islands (MSIs), which are two-dimensional stability islands of the delay-coupling space together with a third dimension (``altitude") encoding for eigenvalues that result in stable AD. We compute the MSFs and show the corresponding MSIs for several common chaotic systems including the R\"{o}ssler, the Lorenz, and Chen's system, and found that it is generally possible to achieve AD and that a nonzero time delay is necessary for the stabilization of the AD states.
\end{abstract}


\maketitle

\section{Introduction}
\textit{Master stability functions} (MSFs), first introduced by Pecora and Carroll \cite{PhysRevLett.80.2109}, provide a framework to evaluate the stability of synchronization of an arbitrary network by mapping its Laplacian eigenvalues to a master stability function which is uniquely determined for a chosen dynamical system (node dynamics), coupling component, and stability measure. This original work has been extended to find the stability of synchronization in networks of coupled dynamical systems with small but arbitrary parametric variations~\cite{jsun,porfiri1}, optimal synchronization in complex networks~\cite{Nishikawa2010PNAS,Ravoori2011PRL}, and synchronization of stochastically coupled chaotic maps~\cite{Juang2008SIAM,porfiri2}. Time delays have important effects on the collective dynamics of coupled oscillators. For networks with time delays, the MSF approach has been adopted to determine the stability of in-phase synchronization and synchronization of specific network topologies of coupled Stuart-Landau oscillators \cite{PhysRevE.81.025205}, synchronization in networks with large coupling delays~\cite{flunkert} and more recently networks with distributed coupling delays~\cite{kyrychko}. In \cite{Gjurchinovski}, the authors use MSFs to investigate the conditions of amplitude death in networks with a different time delay in the coupling versus in the self feedback. 

{\it Amplitude death} (AD) and {\it oscillation death} (OD) are two types of coupling-induced quenching of oscillatory dynamics. In particular, AD typically refers to the stabilization of an otherwise unstable homogeneous (synchronized) fixed point, whereas OD commonly refers to the coupled-induced creation and stabilization of an inhomogeneous (unsynchronized) fixed point~\cite{koseska2013oscillation}. Both AD and OD describe how coupling interactions can quench the otherwise stable oscillations exhibited by uncoupled units. AD and OD have been shown to occur when coupling identical and nonidentical oscillators under various coupling schemes \cite{Erme1990, Miro1990, Redd1998, Pras2005, Koni2003, Karn2007, Pras2010, Resm2011, zou2014emergence}, and have been studied for one-way ring networks \cite{PhysRevE.70.066201}, Erdos-Renyi (ER) random networks \cite{Zou2011}, small-world networks \cite{PhysRevE.68.055103}, and scale-free networks \cite{Liu2009}. AD/OD arise in experimental settings \cite{suresh2014experimental, biwa2015amplitude} and depending on the application, they can be a desirable outcome of the coupled system, such as the stabilization of DC systems \cite{hud_skuf, konishi2015dynamics}, or for circumstances under which AD/OD is undesirable, it has been shown that introducing a proper processing feedback can revive oscillations in coupled nonlinear oscillators thus avoiding the AD/OD regime \cite{zou2013, zou2015restoration, Ghosh2015}. 
An especially insightful result of AD/OD is the emergence of isolated subsets of the two-dimensional delay-coupling strength parameter space ($\tau$ - $\sigma$), which are often called {\it (amplitude/oscillation) death islands} (ADIs/ODIs) \cite{Redd1998, atay2003, zou2012, zou2013, PhysRevE.88.032916}. In particular, for the AD/OD state of a coupled network to be stable it is necessary that the coupling delay $\tau$ and coupling strength $\sigma$ be chosen from within the death islands computed specifically for that network. 

In this paper, we focus on the AD state of networks of delay-coupled oscillators. We adopt the master stability approach to obtain stability regions for AD which can be used for any arbitrary network topology. We introduce the concept of \textit{master stability islands (MSIs)}, which are landscaped stability surfaces obtained from the MSF to include contour (``altitude") information on top of the ADIs. These MSIs can also be viewed as stability slices either in the $\tau$ - $\lambda$ space or in the $\sigma$ - $\lambda$ space, where $\lambda$ a generic parameter associated with the eigenvalues of the coupling matrix $G$. In our numerical experiments we compute MSFs and corresponding stability regions for common chaotic systems such as the R\"{o}ssler system, the Lorenz system, and Chen's system. For all systems, the stabilization of AD requires a nonzero time delay, and the ADIs (and also MSIs) tend to be smaller as the coupling delay increases. For a fixed coupling delay, the set of coupling strengths which correspond to stable AD generally forms a single continuous interval. However, for a fixed coupling strength, the coupling delays that are associated with stable AD form multiple disconnected intervals. Our new MSI computation shows that even within the same ADI, range of stability (visualized as the ``altitude" of different parts of the island) can vary significantly depending on the particular combination of coupling parameters. This highlights the fact that the ADIs alone are not sufficient for the determination of the stability of AD for general networks. Instead, a full MSF (or equivalently, MSI) would be required.

The rest of the paper is organized as follows. In Section \ref{sec:master_stability_approach}, we derive the MSF for  the AD state. We show the conditions needed for diagonalization of the MSF into scalar equations and how the Lambert W-function can be used to find the characteristic roots. Section \ref{sec:def_stability_regions} provides mathematical definitions for the various stability regions. We define and discuss the details of ADIs, MSIs, and stability slices (the slices of the MSIs). In Section \ref{sec:ex_stability_regions}, we show the stability regions in terms of ADIs, MSIs, and stability slices for coupled R\"{o}ssler oscillators, Lorenz oscillators, and Chen's oscillators, respectively, and highlight some common characteristics of the stability regions obtained for such systems. Section \ref{sec:con_and_disc} discusses our results and addresses pertinent issues and observations related to this master stability approach to AD in networks of delay-coupled oscillators.

\section{Master Stability Approach}\label{sec:master_stability_approach}
Consider a network of $n$ identical oscillators with linear delay-coupling. Let $\mathbf{x}_{i}(t)$ and $\mathbf{x}_{i}(t-\tau)$ be the $m$-dimensional vectors of the instantaneous and delayed variables of the $i$th node ($i=1,2,\dots,n$), respectively, and $F(\mathbf{x}_{i})$ the uncoupled dynamics at each node. Then the dynamics of the $i$th node can be written as 
\begin{equation}\label{eq:network}
\dot{\mathbf{x}}_{i} = F(\mathbf{x}_{i}) + \sigma\sum_{j=1}^{n}g_{ij}H[\mathbf{x}_{j}(t-\tau)-\mathbf{x}_{i}(t)],
\end{equation}
where $\sigma$ is the coupling strength, $H\in\mathbb{R}^{m\times m}$ is the node-to-node coupling matrix, and $G=[g_{ij}]_{n\times n}$ is the network coupling matrix where $g_{ij}\geq0$ represents the weight of delayed coupling of node $j$ on node $i$. 
This model has been widely used for the study of coupled oscillator networks with time delays~\cite{flunkert,kyrychko,PhysRevE.81.025205,Dhamala2004,Kinzel2009}.

\subsection{Coupling Matrix}
We make three assumptions about the coupling matrix $G$: (1) the row sum is a constant, that is, $\sum_{j}g_{ij}=c$ for every row $i$; (2) $G$ is diagonalizable; and (3) the eigenvalues of $G$ are all real. These assumptions ensure the existence of a MSF defined on real numbers, which we will discuss later in this section. Without loss of generality, we can absorb the constant row sum $c$ into the coupling strength and simply set $c=1$. One particular example of $G$ that satisfies all the three assumptions is given by $G=D^{-1}A$ where $A$ is any symmetric matrix (e.g., the adjacency matrix of an undirected network with $a_{ij}=1$ if and only if nodes $i$ and $j$ are connected by an edge) and $D$ is a diagonal matrix with diagonal entries defined by $d_{ii}=\sum_{j}a_{ij}$.
The constant row sum of $G$ equaling $1$ poses constraints on the eigenvalues of $G$, which in general can be ordered as:
\begin{equation}
1 = \lambda_{1} \geq \lambda_{2} \geq \dots \geq \lambda_{n} \geq -1,
\end{equation}
where the eigenvalue $\lambda_1=1$ corresponds to the uniform eigenvector $\mathbf{1}=[1,1,\dots,1]^\top$, and the other eigenvalues can be shown to be bounded between $-1$ and $1$ by applying the Ger\v{s}gorin circle theorem~\cite{HornBook} to each row of $G$ together with the non-negativity of the entries of $G$.

\subsection{Stability of Amplitude Death}
We focus on a particular type of synchronization, referred to as {\it amplitude death} (AD), which is characterized by the condition
\begin{equation}\label{eq:od}
\mathbf{x}_1(t)=\mathbf{x}_2(t)=\dots=\mathbf{x}_n(t)=\mathbf{s},
\end{equation}
where $\mathbf{s}$ is an unstable fixed point of the uncoupled system, satisfying $F(\mathbf{s})=0$.
To analyze the stability of the fixed point $\mathbf{s}$, we consider a small arbitrary perturbation $\mathbf{\xi}_{i} =\mathbf{x}_{i}-\mathbf{s}$ of the $i$th node. The time evolution of such a perturbation can be obtained by linearization of Eq.~\eqref{eq:network}, giving rise to a variational equation
\begin{equation}\label{eq:variational1}
  \dot{\mathbf{\xi}}_{i} = DF(\mathbf{s})\mathbf{\xi}_{i}(t) + \sigma\sum_{j=1}^{n}g_{ij}H[\mathbf{\xi}_{j}(t-\tau) - 
  \mathbf{\xi}_{i}(t)].
  \end{equation}
  Let $\mathbf{\xi} = (\mathbf{\xi}_{1}, \mathbf{\xi}_{2},\dots ,
\mathbf{\xi}_{n})$ be the collection of node variations for all the nodes in the network. Then  Eq. (\ref{eq:variational1}) can be expressed in matrix form as
\begin{equation}\label{eq:matrix_form}
  \dot{\mathbf{\xi}} = [I_{n} \otimes (D
  F(\mathbf{s})-\sigma H)]  \mathbf{\xi}(t)  +
  \sigma (G \otimes 
  H) \mathbf{\xi} (t-\tau),
  \end{equation}
where $\otimes$ represents the Kronecker product. The goal is to  diagonalize Eq. (\ref{eq:matrix_form}) so that the stability of each mode can be analyzed separately. 
In order to accomplish this, $G$ must be diagonalized, as
\begin{equation}
\Lambda = P^{-1}GP =
\text{diag}(\lambda_{1},\lambda_{2},\dots,\lambda_{n}),
\end{equation}
where the $\lambda_{i}$'s are the eigenvalues of $G$.
By applying the change of variables $\mathbf{\eta} = (P^{-1}
\otimes I_{m}) \mathbf{\xi}$, Eq. (\ref{eq:matrix_form}) becomes
\[
  \dot{\mathbf{\eta}} = [I_{N} \otimes (D
  F(\mathbf{s})-\sigma H)]\mathbf{\eta}(t) +
  \sigma (\Lambda \otimes H) \mathbf{\eta}(t - \tau).
  \]
Thus, there are $n$ independent modes of the form
\begin{equation}\label{eq:etai}
  \dot{\mathbf{\eta}}_{i} = (DF(\mathbf{s}) - \sigma
  H) \mathbf{\eta}_{i}(t) + \sigma
  \lambda_{i} H \mathbf{\eta}_{i}(t-\tau ).
\end{equation}

\subsection{Master Stability Functions}
Since the form of the Eq.~\eqref{eq:etai} remains the same for each block, we can define a master stability equation as
  \begin{equation}\label{eq:matrix_msf}
\dot{\mathbf{\zeta}} = (DF(\mathbf{s}) -
  \sigma H)\mathbf{\zeta}(t) + \sigma\lambda H \mathbf{\zeta}(t-\tau).
  \end{equation}
  
For a given system specified by $F$ and $H$ and the AD state $\mathbf{s}$, we will denote the maximum real part of the characteristic roots associated with Eq.~\eqref{eq:matrix_msf} by $\Omega(\tau,\sigma,\lambda)$, which depends on three inputs: $\tau$ (coupling delay), $\sigma$ (coupling strength), and $\lambda$ (a generic parameter associated with the eigenvalues of $G$). We refer to this function $\Omega:\mathbb{R}\times\mathbb{R}\times\mathbb{R}\rightarrow\mathbb{R}$ as the {\it master stability function} (MSF) for the AD state $\mathbf{s}$. In particular, for an arbitrary matrix $G$ which satisfies the three assumptions outlined at the beginning of this section and whose eigenvalues are $\{\lambda_i\}_{i=1}^{n}$, a sufficient condition for the AD state to be stable is given by
\begin{equation}
	\max_{1\leq i\leq n}\Omega(\tau,\sigma,\lambda_i)<0,
\end{equation}
which can be solely determined from the MSF.
  
Below we show how to compute the characteristic roots of the multivariate differential Eq.~\eqref{eq:matrix_msf} by further decomposing it into scalar differential equations and utilizing the Lambert W-function. Such decomposition requires 
the matrices $DF(\mathbf{s})$ and $H$ to commute and be both diagonalizable. Since commuting matrices have the same set of eigenvectors, there exists an invertible matrix $Q$ that simultaneously diagonalize both $DF(\mathbf{s})$ and $H$, as
\begin{equation}\label{eq:DFHdiag}
\begin{cases}
    DF(s) = QM^{(DF)}Q^{-1},\\
    H = Q M^{(H)}Q^{-1},
\end{cases}
\end{equation}
where $M^{(DF)}$ and $M^{(H)}$ are diagonal matrices whose diagonal elements $\{\mu^{(DF)}_\ell\}_{\ell=1}^{m}$ and $\{\mu^{(H)}_\ell\}_{\ell=1}^{m}$ are the set of eigenvalues of $DF(\mathbf{s})$ and $H$, respectively. The diagonalizations in Eq.~\eqref{eq:DFHdiag} can be used along with the change of variable $\mathbf{\psi}=Q^{-1}\mathbf{\zeta}$ to transform Eq. (\ref{eq:matrix_msf}) into $m$ decoupled scalar equations, for $\ell=1,2,\dots,m$:
\begin{equation}\label{eq:scalar_mse}
\dot{\psi}_{\ell} = \left(\mu^{(DF)}_{\ell} -\sigma\mu^{(H)}_\ell\right)\psi(t) + \sigma\lambda\mu^{(H)}_\ell\psi(t-\tau).
\end{equation}
This scalar master stability equation has the corresponding characteristic equation
\begin{equation}\label{eq:scalar_char}
\mu = \mu^{(DF)}_{\ell} -\sigma\mu^{(H)}_\ell + \sigma\lambda\mu^{(H)}_\ell e^{-\mu \tau}.
\end{equation}
The solution of Eq.~\eqref{eq:scalar_mse} satisfies $\psi_\ell(t)\rightarrow0$ as $t\rightarrow\infty$ if all the real parts of the roots of Eq. (\ref{eq:scalar_char}) are negative. When this occurs, the AD state of the coupled network system becomes stable.
The roots of Eq. (\ref{eq:scalar_char}) can be expressed as
\begin{equation}\label{eq:char}
\mu = \mu_\ell^{(DF)}-\sigma\mu^{(H)}_\ell + \frac{1}{\tau}W(\sigma \lambda \tau e^{-\tau(\mu_\ell^{(DF)}-\sigma \mu^{(H)}_\ell)}),
\end{equation}
where $W(\cdot)$ denotes the Lamber W function~\cite{878632}, which is in fact a multivalued inverse of the function $w\mapsto we^w$~\cite{Corless1996}.
While the Lambert W-function has an infinite number of branches, Shinozaki and Mori proved that the principle branch always determines the stability of a scalar linear delay differential equation~\cite{Shinozaki20061791}.
Thus, the MSF value $\Omega(\tau,\sigma,\lambda)$ can be computed as the maximum real part of the root of Eq.~\eqref{eq:char} using the principle branch of the Lambert W-function, maximized over the indices $\ell=1,2,\dots,m$ after diagonalization of the matrices $DF(\mathbf{s})$ and $H$.

In the case where the matrices $DF(\mathbf{s})$ and $H$ do not commute or are not both diagonalizable, the matrix master stability Eq.~\eqref{eq:matrix_msf} cannot be diagonalized. Computation of the characteristic roots in this case is more involved. As noted in~\cite{MichielsBook}, the problem can be formulated as a nonlinear eigenvalue and numerically solved using spectral discretization. Software packages such as the DDE-BIFTOOL~\cite{DDE_BIFTOOL1, DDE_BIFTOOL2} compute these roots via a linear multi-step method as detailed in the manual~\cite{DDE_BIFTOOL2}.

\section{Definition of Stability Regions}\label{sec:def_stability_regions}
\subsection{Amplitude Death Islands}
In the $\tau$ - $\sigma$ space, note that the condition of constant row sum equalling one of the matrix $G$ implies that $\lambda=1$ is always an eigenvalue. Therefore, a {\it necessary} condition for system~\eqref{eq:network} to have stable AD is
\begin{equation}\label{eq:odi}
\Omega(\tau,\sigma,\lambda=1)<0.
\end{equation}
The set of parameter combinations $(\tau,\sigma)$ for a region in the $\tau$ - $\sigma$ space, which is typically made up of isolated regions visually looking like islands (see Figs.~\ref{fig:rossler_od},~\ref{fig:lorenz_od}, and~\ref{fig:chen_od} for amplitude death islands of the R\"{o}ssler system, Lorenz system, and Chen's system, respectively.) For this reason, these island-like regions are often called amplitude death islands (ADIs) (or amplitude death islands) in the literature \cite{Redd1998, atay2003, zou2012, zou2013, PhysRevE.88.032916}. Note, however, that in some papers the ADIs are defined for specific types of networks as the parameter combinations of $\tau$ and $\sigma$ under which the particular networks under consideration have stable AD \cite{PhysRevE.84.066208, zou2009partial, zou2010eliminating, zou2012, zou2013}.

\subsection{Master Stability Islands}
The way to interpret the ADI is that it imposes a necessary condition for a network to have stable AD by requiring the delay and coupling strength to be chosen from one of these islands. However, for a given network the ADIs alone are not enough/sufficient to determine the stability of AD. The reason is that there might exist eigenvalues $\lambda\neq1$ of $G$ which correspond to positive values (unstable regimes) of the MSF. To account for the influence of the eigenvalues in addition to the parameters $\tau$ and $\sigma$, we define, for each parameter combination $(\tau,\sigma)$ inside an ADI, a unique set of $\lambda$'s for which the MSF is negative (stable). Such a set is given by
\begin{equation}\label{eq:msi}
I_\lambda(\sigma,\tau)=\{\lambda|\Omega(\sigma,\tau,\lambda)<0\}\cap[-1,1].
\end{equation}
From all numerical experiments that we have performed, the stability set $I_\lambda$ always takes the form of a continuous interval $[a,1]$ (although we were not able to prove this). Thus, for each ADI there is a ``landscape" defined by using the length of the stability intervals $I_\lambda$ as ``altitudes" to capture the range of eigenvalues that fall within the stability region of AD. This renewed concept of stability islands gives rise to what we call {\it master stability islands} (MSIs), as the MSIs (just like the MSFs) suffice to determine the stability of AD of an arbitrary network. Examples of MSIs are shown in Figs.~\ref{fig:rossler_msi}, \ref{fig:lorenz_msi1}, and \ref{fig:chen_msi}, with detailed descriptions provided in Section \ref{sec:ex_stability_regions}.

\subsection{Stability Slices}
Another way to visualize the stability regions are to plot them in the $\sigma$ - $\lambda$ space upon different choices of the delay parameter $\tau$ (as shown in Figs.~\ref{fig:rossler_sigma},~\ref{fig:lorenz_k}, and~\ref{fig:chen_sigma}), or in the $\tau$ - $\lambda$ space upon different choices of the coupling strength $\sigma$ (as shown in Figs.~\ref{fig:rossler_tau},~\ref{fig:lorenz_tau}, and~\ref{fig:chen_tau}). These regions can be thought of as slices of the master stability surface living in $\tau$ - $\sigma$ - $\lambda$ space.
The stability regions where the MSF values are negative are shaded in gray in these figures.
Recall that the eigenvalues of $G$ must lie between $-1$ and $1$. Thus, in these figures the dashed horizontal lines at $\lambda=\pm 1$ are used to indicate the guaranteed maximum value and possible minimum value of $\lambda$, giving rise to subregions marked by the slanted lines. The rest of the gray regions are simply not realizable for the model that we consider.

\section{Examples of Stability Regions}\label{sec:ex_stability_regions}
In this section, we plot the numerically determined stability regions in the forms of stability islands and stability slices for the R\"{o}ssler system, the Lorenz system, and Chen's system in the various parameter spaces $\tau$ - $\sigma$, $\tau$ - $\lambda$, and $\sigma$ - $\lambda$. 

\subsection{R\"{o}ssler System}
The R\"{o}ssler system \cite{Rossler1976397} is given by
\begin{equation}\label{eq:rossler}
\begin{cases}
\dot{x} & =  -y-z,  \\
\dot{y} & =  x+ay, \\
\dot{z} & =  b+(x-c)z,
\end{cases}
\end{equation}
and has two real fixed points under the condition of $c^2>4ab$:
\begin{equation}
  \begin{aligned}
    x^{*}_{1,2} &= \frac{c \pm \sqrt{c^{2}-4ab}}{2a}\\
    y^{*}_{1,2} &= \frac{-c \pm \sqrt{c^{2}-4ab}}{2}\\
    z^{*}_{1,2} &= \frac{c \pm \sqrt{c^{2}-4ab}}{2}.
    \end{aligned}
\end{equation}
The Jacobian matrix is
\begin{equation}
  DF = 
  \begin{bmatrix}
    0 & -1 & -1 \\
    1 & a & 0 \\
    z & 0 & x-c
    \end{bmatrix}.
\end{equation}
For the parameter values $a=0.15$, $b=0.2$, and $c=10$, one of the fixed points is at $x^{*} = 0.003$, $y^{*} = -0.02$, and $z^{*} = 0.02$. The eigenvalues of the Jacobian evaluated at this fixed point are $\mu_{1} = 0.0740 + 0.9972i$, $\mu_{2} = 0.0740 - 0.9972i$, and $\mu_{3} = -9.9950$.

\begin{figure}[h!] 
\includegraphics[scale=1]{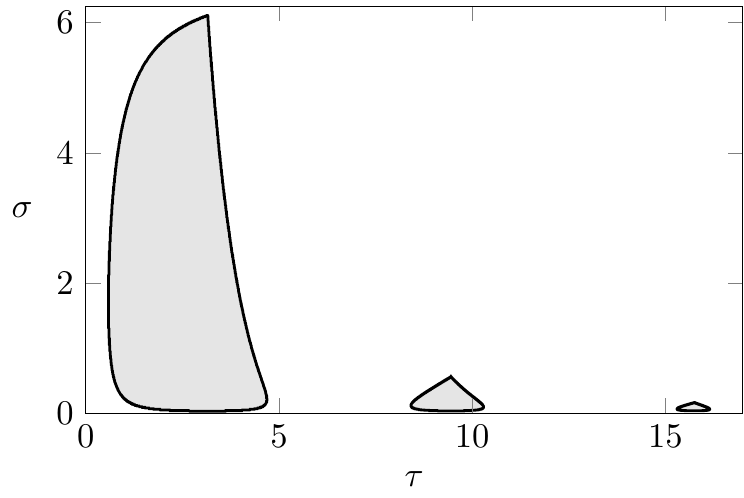} 
\caption{\label{fig:rossler_od} R\"{o}ssler system amplitude death islands (ADIs). The system is defined by Eq.~\eqref{eq:rossler} with parameter values $a=0.15$, $b=0.2$, and $c=10$ and $x \rightarrow x$, $y \rightarrow y$, and $z \rightarrow z$ coupling. The ADIs, given by Eq.~\eqref{eq:odi}, are obtained from the numerically determined master stability function (MSF). }
\end{figure}

\begin{figure}[h!] 
\includegraphics[scale=1]{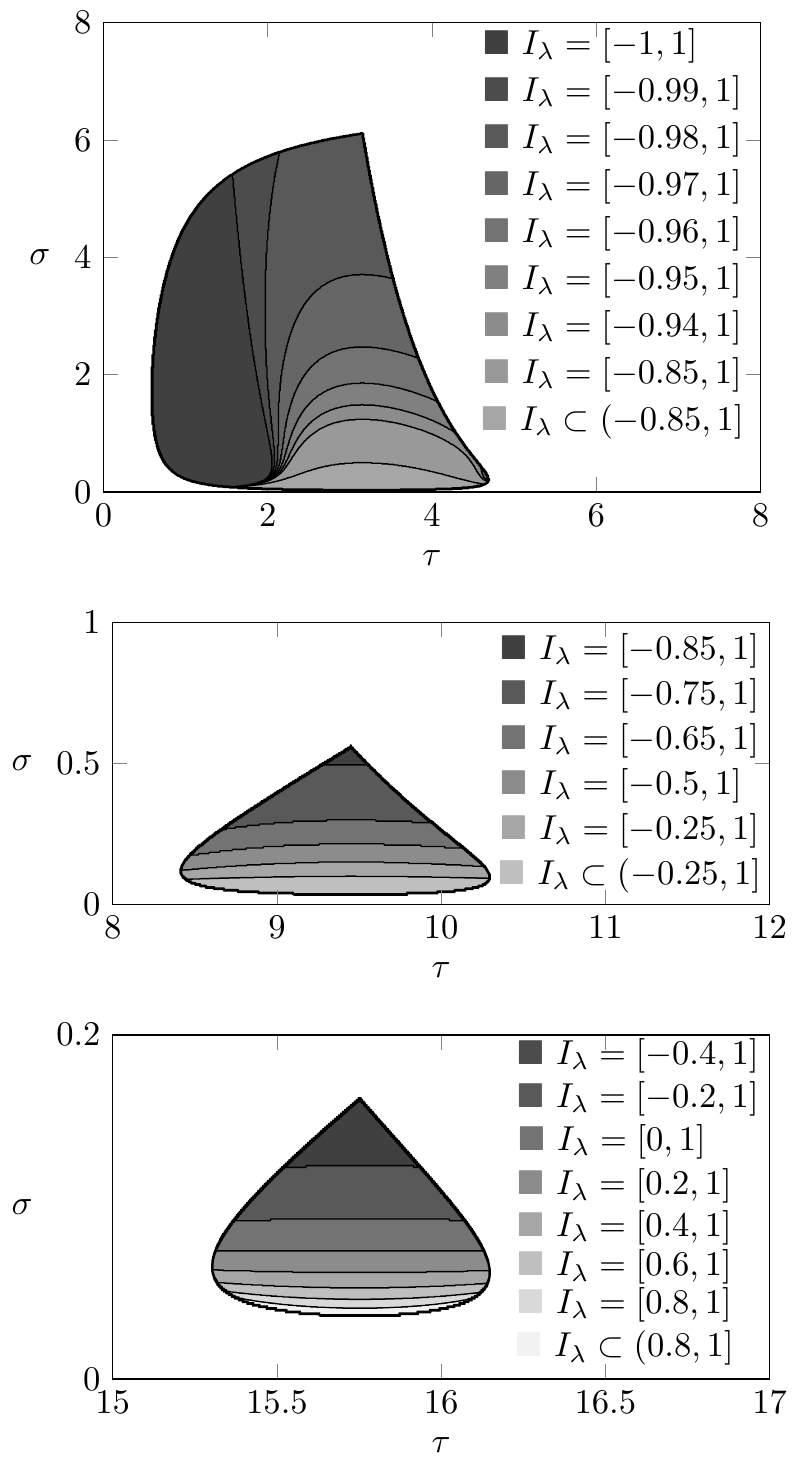} 
\caption{\label{fig:rossler_msi} R\"{o}ssler system master stability islands (MSIs). The system is defined by Eq.~\eqref{eq:rossler} with parameter values $a=0.15$, $b=0.2$, and $c=10$ and $x \rightarrow x$, $y \rightarrow y$, and $z \rightarrow z$ coupling. The MSIs, given by Eq.~\eqref{eq:msi}, are obtained from the numerically determined master stability function (MSF).}
\end{figure}

\begin{figure}[h!]
\hspace{-3mm} \includegraphics[scale=1]{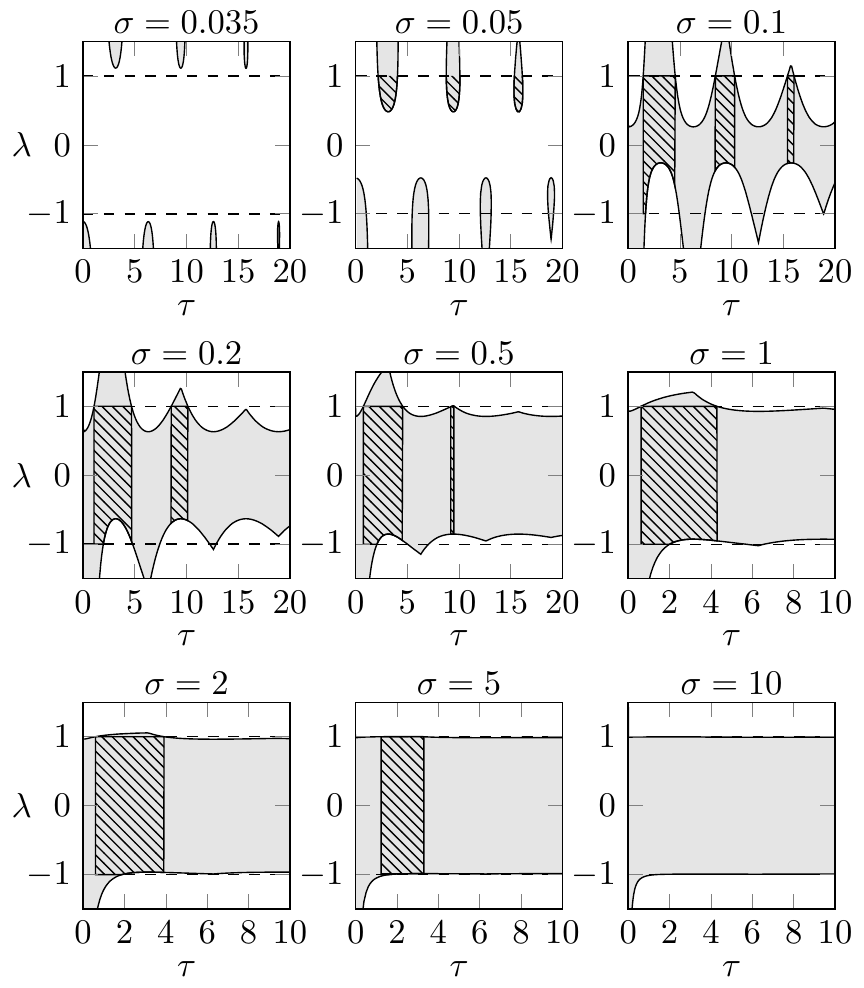} 
\caption{\label{fig:rossler_sigma} R\"{o}ssler system stability slices in the $\tau$ - $\lambda$ space. The system is defined by Eq.~\eqref{eq:rossler} with parameter values $a=0.15$, $b=0.2$, and $c=10$ and $x \rightarrow x$, $y \rightarrow y$, and $z \rightarrow z$ coupling. In each panel, the active stability regions (regions that contain $\lambda=1$) are filled with the slanted lines.}
\end{figure}

\begin{figure}[h!]
\hspace{-3mm} \includegraphics[scale=1]{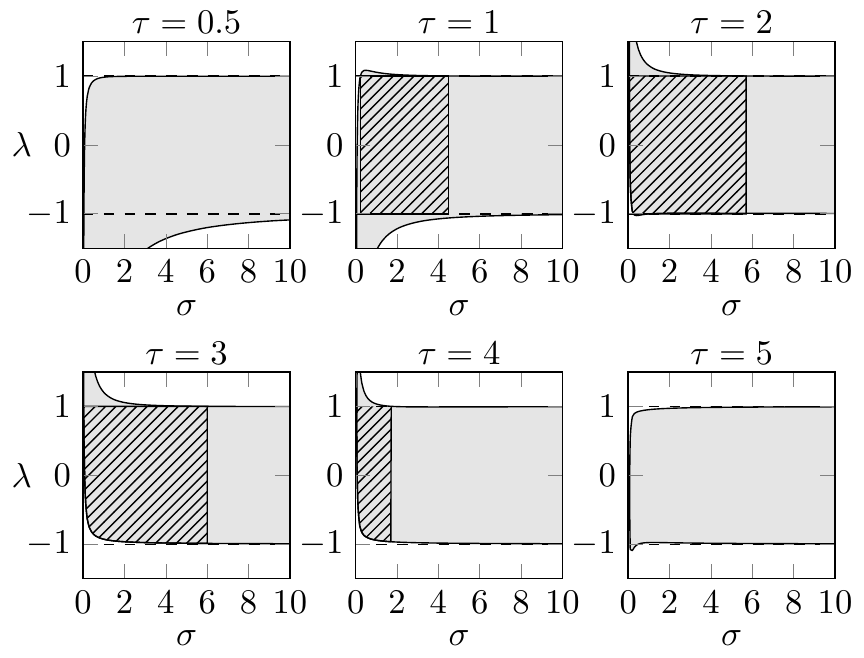}
  \caption{\label{fig:rossler_tau}  R\"{o}ssler system stability slices in the $\sigma$ - $\lambda$ space. The system is defined by Eq.~\eqref{eq:rossler} with parameter values $a=0.15$, $b=0.2$, and $c=10$ and $x \rightarrow x$, $y \rightarrow y$, and $z \rightarrow z$ coupling. In each panel, the active stability regions (regions that contain $\lambda=1$) are filled with the slanted lines.}
\end{figure}

Under the above-mentioned parameters and coupling in all variables, that is $H=I$, the R\"{o}ssler system has three amplitude death islands over the range of $(\tau,\sigma)\in[0, 20]\times[0, 20]$. These ADIs, which are obtained from the MSF according to Eq.~\eqref{eq:odi}, are shown in Fig.~\ref{fig:rossler_od}. Note that none of the ADIs touch the $\tau=0$ axis, indicating that a nonzero time delay is necessary for the stabilization of AD.
The corresponding MSIs whose ``atitudes" are determined by Eq.~\eqref{eq:msi} are shown in Fig.~\ref{fig:rossler_msi}. Here we found that the size of the islands becomes smaller as the time delay $\tau$ increases, a phenomenon also observed for the Lorenz oscillators which will be presented later. Fig. \ref{fig:rossler_sigma} shows the stability slices of the coupled system for fixed coupling strengths, and Fig.~\ref{fig:rossler_tau} shows the stability slices of the coupled system for fixed delays. For the system that we consider, the coupling matrix $G$ has constant row sums and consequently $\lambda=1$ is always an eigenvalue. Thus, for a region in the stability slices shown in Figs.~\ref{fig:rossler_sigma} and~\ref{fig:rossler_tau} to be active in determining the stability of AD, it has to include $\lambda=1$. In the figures, we highlight these active stability regions by filling them in with slanted lines. For fixed coupling strength $\sigma$, the size and number of active stability regions in the $\tau$ - $\lambda$ space are found to depend (non-monotonically) on the value of the coupling strength as shown in Fig.~\ref{fig:rossler_sigma}. On the other hand, for fixed delay $\tau$, there is generally either none or a single connected stability region in the $\tau$ - $\lambda$ space as shown in Fig.~\ref{fig:rossler_tau}.

\subsection{Lorenz System}
The Lorenz system \cite{lorenz1963} is given by
\begin{equation}\label{eq:lorenz}
\begin{cases}
\dot{x} & =  a(y-x)  \\
\dot{y} & =  x(r -z)-y \\
\dot{z} & =  xy- b z.
\end{cases}
\end{equation}
If $r < 1$, then the origin is the only fixed point. For $r > 1$ there exists two fixed points 
\begin{equation}
  \begin{aligned}
    x^{*}_{1,2} &= \pm \sqrt{\b(r-1)}\\
    y^{*}_{1,2} &= \pm \sqrt{\b(r-1)}\\
    z^{*}_{1,2} &= r-1.
    \end{aligned}
\end{equation}
The Jacobian matrix is
\begin{equation}
  DF = 
  \begin{bmatrix}
    -a & a & 0 \\
    r -z & -1 & -x \\
    y & x & -b
    \end{bmatrix}.
\end{equation}
For the parameter values $a=10$, $r=8/3$, and $b=28$, one of the fixed points is at $x^{*} = 0.485$, $y^{*} = 0.485$, and $z^{*} = 27$. The eigenvalues of the Jacobian are $\mu_{1} = 0.0939 + 10.1945i$, $\mu_{2} = 0.0939 - 10.1945i$, and $\mu_{3} = -13.8546$.

\begin{figure}[h!] 
\includegraphics[scale=1]{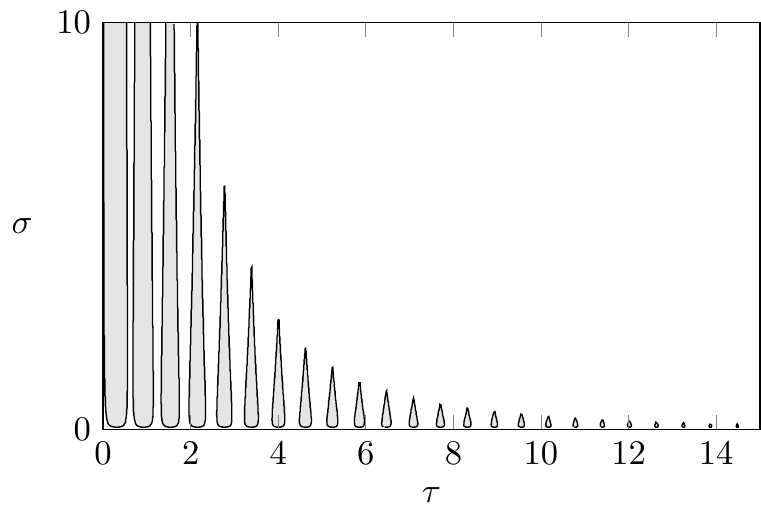} 
\caption{\label{fig:lorenz_od} Lorenz system amplitude death islands (ADIs). The system is defined by Eq.~\eqref{eq:lorenz} with parameter values $a=10$, $r=8/3$, and $b=28$ and $x \rightarrow x$, $y \rightarrow y$, and $z \rightarrow z$ coupling. The ADIs, given by Eq.~\eqref{eq:odi}, are obtained from the numerically determined master stability function (MSF).}
\end{figure}

\begin{figure}[h!] 
\includegraphics[scale=1]{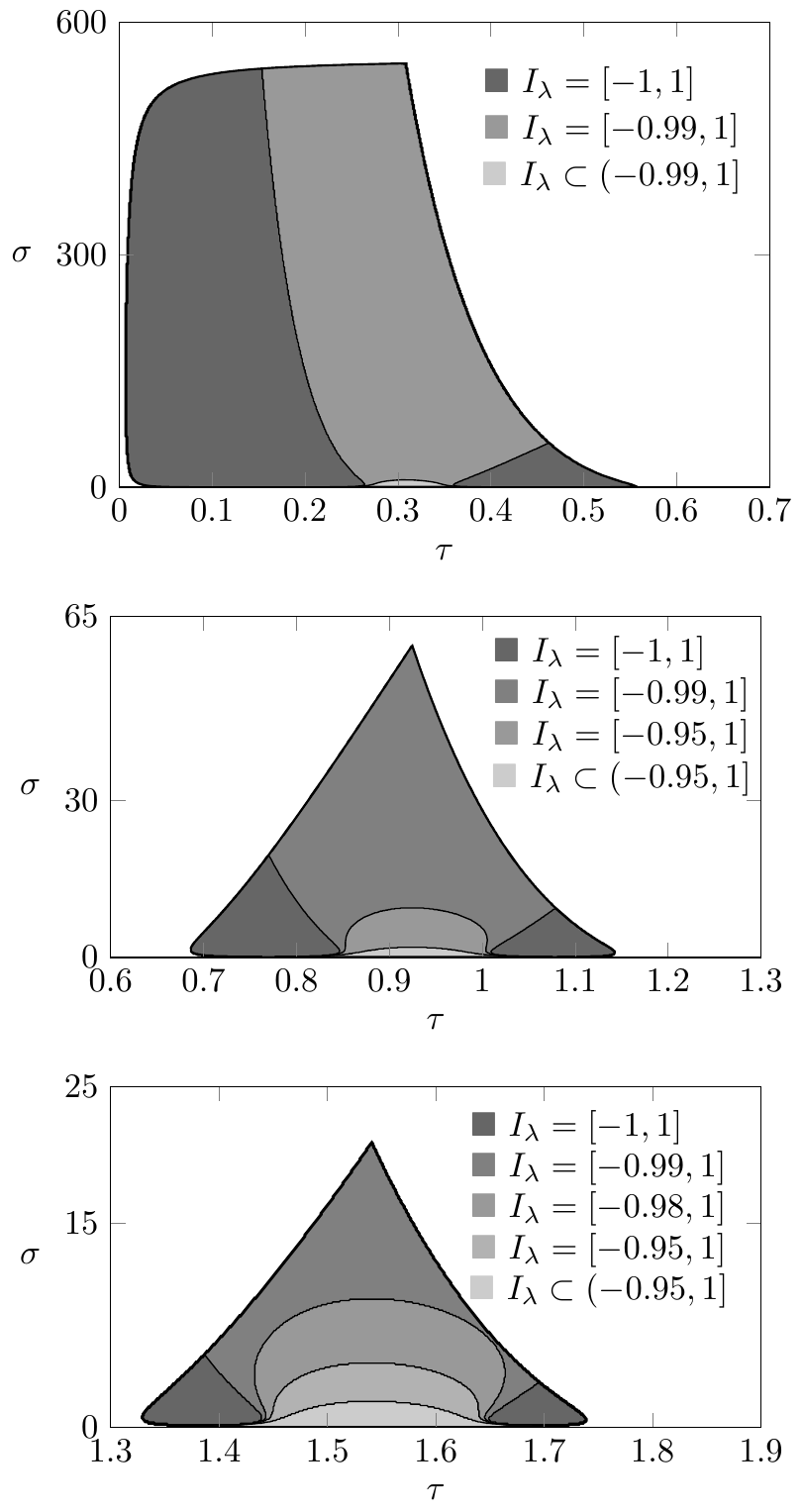} 
\caption{\label{fig:lorenz_msi1} Lorenz system first, second, and third master stability islands (MSIs). The system is defined by Eq.~\eqref{eq:lorenz} with parameter values $a=10$, $r=8/3$, and $b=28$ and $x \rightarrow x$, $y \rightarrow y$, and $z \rightarrow z$ coupling. The MSIs, given by Eq.~\eqref{eq:msi}, are obtained from the numerically determined master stability function (MSF).}
\end{figure}

\begin{figure}[h!]
\hspace{-5mm}   \includegraphics[scale=1]{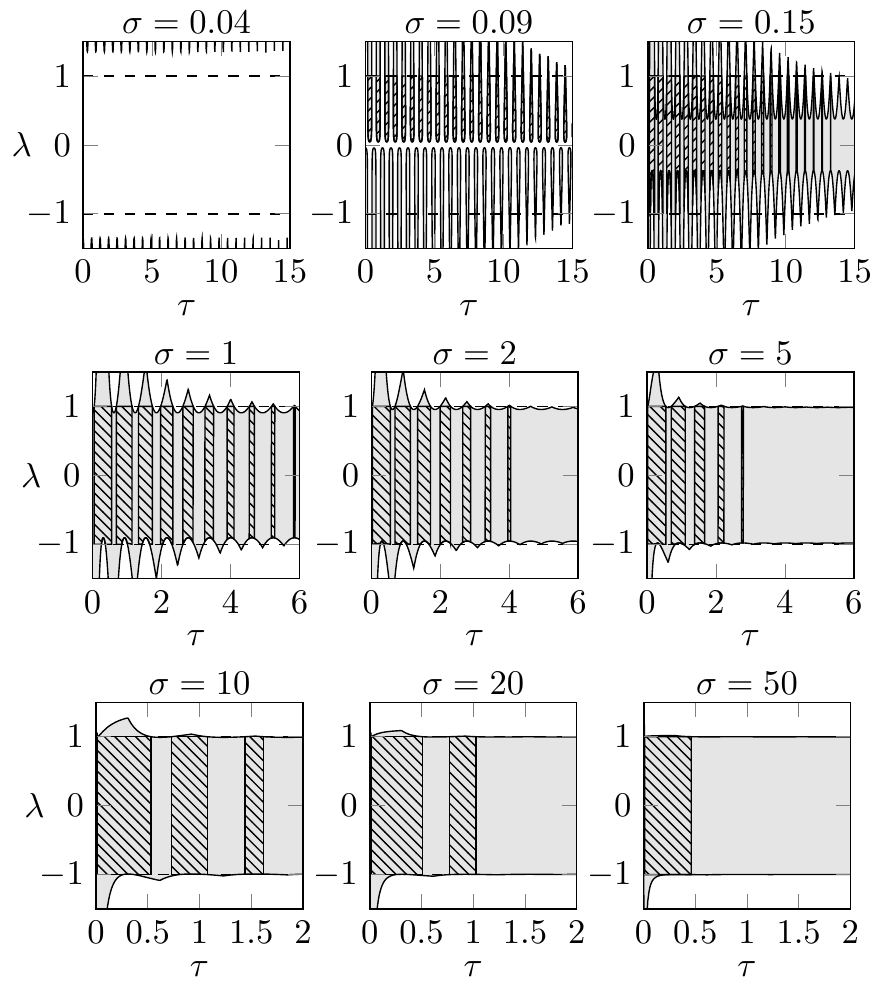}
  \caption{\label{fig:lorenz_k} Lorenz system stability slices in the $\tau$ - $\lambda$ space. The system is defined by Eq.~\eqref{eq:lorenz} with parameter values $a=10$, $r=8/3$, and $b=28$ and $x \rightarrow x$, $y \rightarrow y$, and $z \rightarrow z$ coupling. In each panel, the active stability regions (regions that contain $\lambda=1$) are filled with the slanted lines.}
\end{figure}

\begin{figure}[h!]
\hspace{-5mm}   \includegraphics[scale=1]{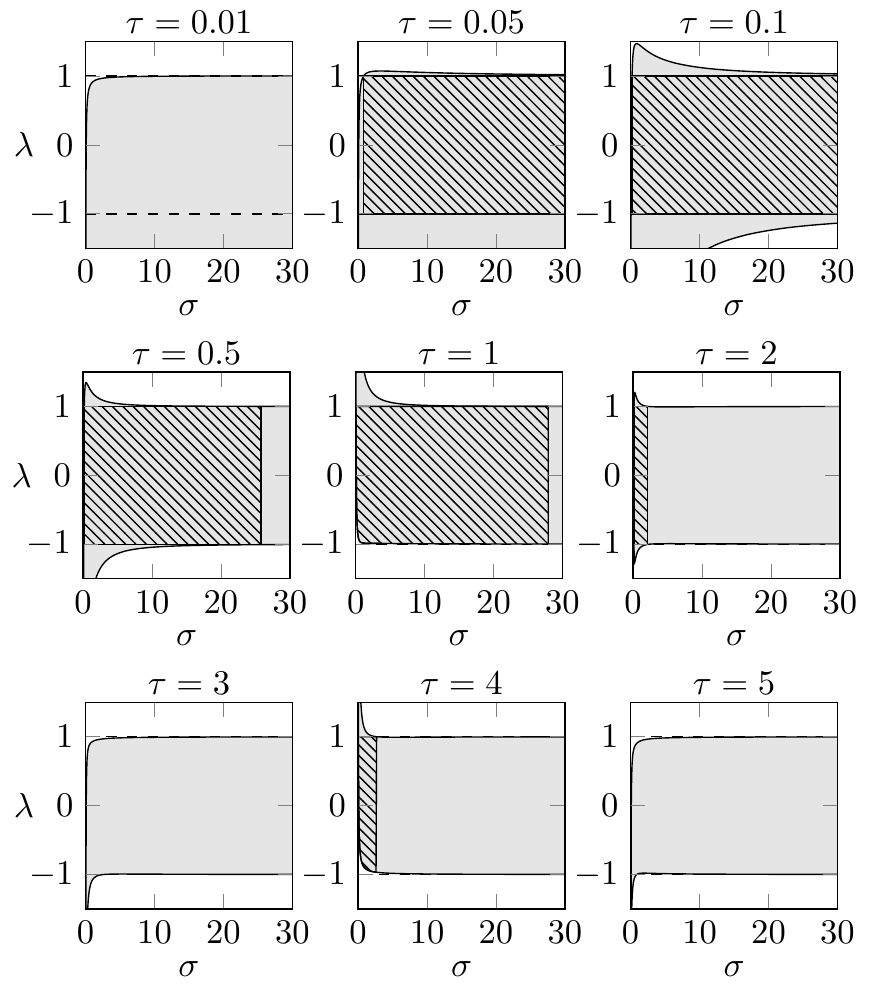}
  \caption{\label{fig:lorenz_tau} Lorenz system stability slices in the $\sigma$ - $\lambda$ space. The system is defined by Eq.~\eqref{eq:lorenz} with parameter values $a=10$, $r=8/3$, and $b=28$ and $x \rightarrow x$, $y \rightarrow y$, and $z \rightarrow z$ coupling. In each panel, the active stability regions (regions that contain $\lambda=1$) are filled with the slanted lines.}
\end{figure}

Under the above-mentioned parameters and coupling in all variables, that is $H=I$, the Lorenz system has twenty-four amplitude death islands over the range of $(\tau,\sigma)\in[0, 15]\times[0, 600]$. These ADIs are shown in Fig.~\ref{fig:lorenz_od} and the first three corresponding MSIs are shown in Fig.~\ref{fig:lorenz_msi1}. Similar to the R\"{o}ssler system, here all the ADIs are away from $\tau=0$, implying the necessity of having a nonzero time delay in order for the AD state to be stable.
As mentioned above, we found that the size of these islands becomes smaller as the time delay $\tau$ is increased.  Fig.~\ref{fig:lorenz_k} shows the stability slices of the coupled system for fixed coupling strengths, and Fig.~\ref{fig:lorenz_tau} shows the stability slices of the coupled system for fixed delays. Corresponding to the results with the R\"{o}ssler system, for fixed coupling strength $\sigma$, the size and number of active stability regions in the $\tau$ - $\lambda$ space are found to depend (non-monotonically) on the value of the coupling strength as shown in Fig.~\ref{fig:lorenz_k}, and for fixed delay $\tau$, there is generally either none or a single connected stability region in the $\tau$ - $\lambda$ space as shown in Fig.~\ref{fig:lorenz_tau}. As with the R\"{o}ssler system figures, the regions with the slanted lines represent the active stability regions. 

\subsection{Chen's System}
Chen's system \cite{chen1999yet} is given by
\begin{equation} \label{eq:chen}
\begin{cases}
\dot{x} & =  a(y-x)  \\
\dot{y} & =  (c-a-z)x+cy\\
\dot{z} & =  xy-\beta z.
\end{cases}
\end{equation}
Chen's system always has the fixed point $(0,0,0)$. If $\beta(2c-a)>0$, then there exists two more fixed
points
\begin{equation}
  \begin{aligned}
    x^{*}_{1,2} &= \pm \sqrt{\beta(2c-a)}\\
    y^{*}_{1,2} &= \pm \sqrt{\beta(2c-a)}\\
    z^{*}_{1,2} &= 2c-a.
    \end{aligned}
\end{equation}
The Jacobian matrix is
\begin{equation}
  DF = 
  \begin{bmatrix}
    -a & a & 0 \\
    c-a-z & c & -x \\
    y & x & -\beta
    \end{bmatrix}.
\end{equation}
For the parameter values $a=35$, $c=28$, and $\beta=8/3$, one of the fixed points is at $x^{*} = 7.483$, $y^{*} = 7.483$, and $z^{*} = 21$. The eigenvalues of the Jacobian are $\mu_{1} = 4.0769 + 14.2601i$, $\mu_{2} = 4.0769 - 14.2601i$, and $\mu_{3} = -17.8205$.

\begin{figure}[h!] 
\includegraphics[scale=1]{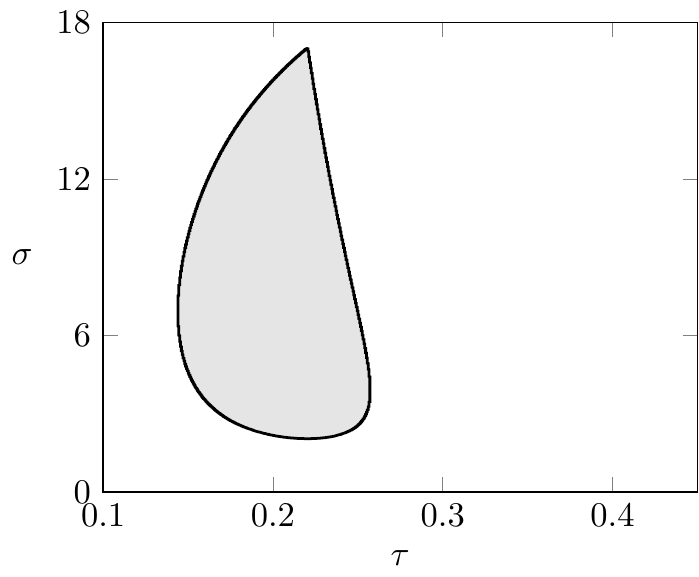} 
\caption{\label{fig:chen_od} Chen's System amplitude death island (ADI). The system is defined by Eq.~\eqref{eq:chen} with parameter values $a=35$, $c=28$, and $\beta=8/3$ and $x \rightarrow x$, $y \rightarrow y$, and $z \rightarrow z$ coupling. The ADI, given by Eq.~\eqref{eq:odi}, is obtained from the numerically determined master stability function (MSF).}
\end{figure}

\begin{figure}[h!] 
\includegraphics[scale=1]{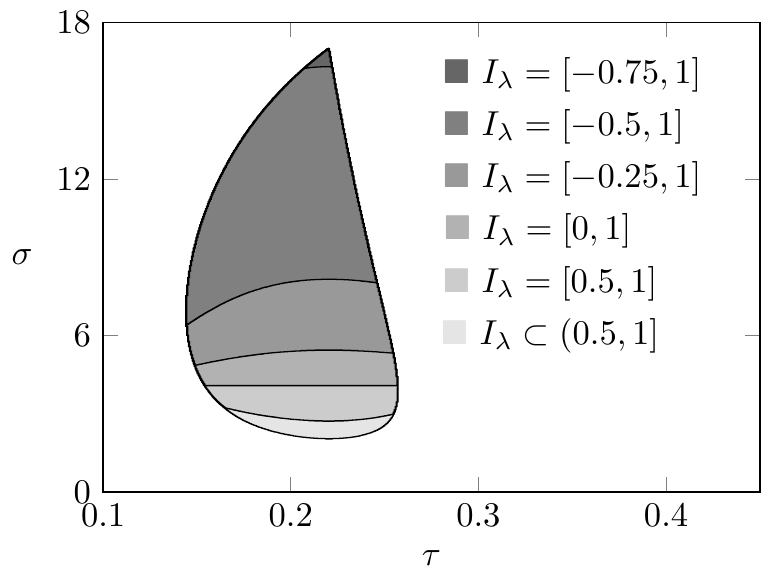} 
\caption{\label{fig:chen_msi} Chen's system master stability island MSI. The system is defined by Eq.~\eqref{eq:chen} with parameter values $a=35$, $c=28$, and $\beta=8/3$ and $x \rightarrow x$, $y \rightarrow y$, and $z \rightarrow z$ coupling. The MSI, given by Eq.~\eqref{eq:msi}, is obtained from the numerically determined master stability function (MSF).}
\end{figure}

\begin{figure}[h!]
\hspace{-5mm}   \includegraphics[scale=1]{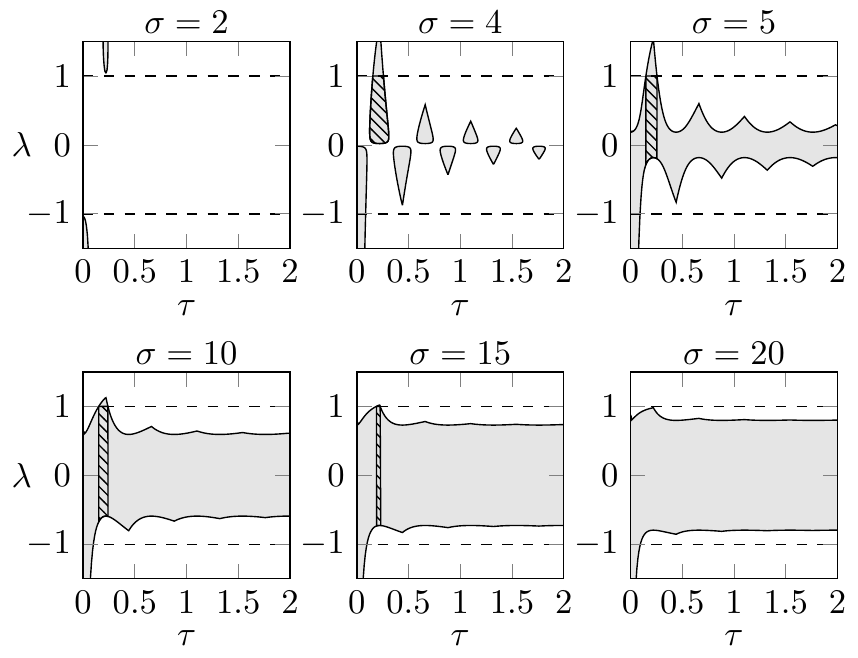}
\caption{\label{fig:chen_sigma} Chen's system stability slices in the $\tau$ - $\lambda$ space. The system is defined by Eq.~\eqref{eq:chen} with parameter values $a=35$, $c=28$, and $\beta=8/3$ and $x \rightarrow x$, $y \rightarrow y$, and $z \rightarrow z$ coupling. In each panel, the active stability regions (regions that contain $\lambda=1$) are filled with the slanted lines.}
\end{figure}

\begin{figure}[h!]
\hspace{-4mm} \includegraphics[scale=1]{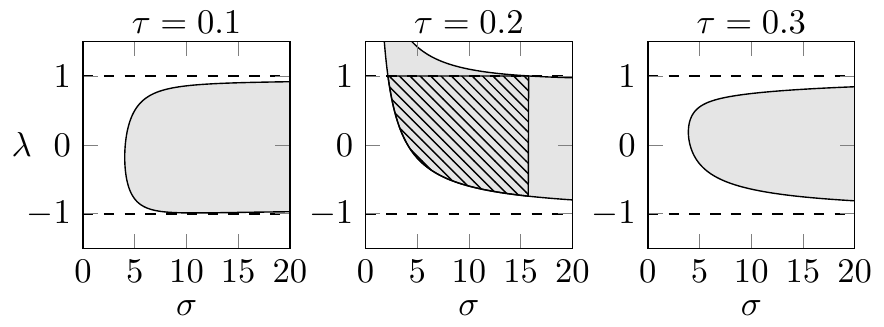}
\caption{\label{fig:chen_tau} Chen's system stability slices in the $\sigma$ - $\lambda$ space. The system is defined by Eq.~\eqref{eq:chen} with parameter values $a=35$, $c=28$, and $\beta=8/3$ and $x \rightarrow x$, $y \rightarrow y$, and $z \rightarrow z$ coupling. In each panel, the active stability regions (regions that contain $\lambda=1$) are filled with the slanted lines.}
\end{figure}

Under the above-mentioned parameters and coupling in all variables, that is $H=I$, Chen's system has one amplitude death island over the range of $(\tau,\sigma)\in[0, 15]\times[0, 30]$. This ADI is shown in Fig.~\ref{fig:chen_od} with its corresponding MSI shown in Fig.~\ref{fig:chen_msi}. Similar to the R\"{o}ssler system and the Lorenz system, the ADI for the Chen's system only exists for $\tau\neq0$, which suggests that a nonzero time delay is required for the stabilization of the AD state. One interesting property that is observed for the Chen's system (but not the R\"{o}ssler or the Lorenz system), is the fact that the MSI of the Chen's system does not contain any region corresponding to $I_{\lambda} = [-1,1]$.  Therefore, particular networks, such as ring networks with an even number of nodes, will not display stable oscillation with the typical setup of $G_{ij}=A_{ij}/k_i$, where $A=[A_{ij}]$ is the adjacency matrix of the (undirected) network and $k_i=\sum_{j}A_{ij}$ is the in-degree of node $i$. This is because under this setting, $\lambda=-1$ is always an eigenvalue of $G$ for these networks \cite{Zou2011}. Fig.~\ref{fig:chen_sigma} shows the stability slices of the coupled system for fixed coupling strengths, and Fig.~\ref{fig:chen_tau} shows the stability slices of the coupled system for fixed delays. As with the figures above, the regions with the slanted lines represent the active stability regions. 

\section{Conclusions and Discussion}\label{sec:con_and_disc}
In this paper we derived the MSF for the AD state using the approach of Pecora and Carroll  \cite{PhysRevLett.80.2109} slightly modified for delay differential equations. We provided the conditions under which the vector characteristic equations could be decoupled into scalar ones and showed how to use the Lambert W-function to calculate their roots in order to obtain the corresponding MSF. 
We then introduced the concept of MSIs, which are a visual representation of the ADIs together with ``altitude" information that can be used to determine the stability of AD of arbitrary networks. Next we defined stability slices (slices of the MSI) and demonstrated how they can be an insightful way to view the corresponding parameter spaces. Finally, we provided extensive numerical experiments from which MSFs and corresponding stability regions including MSIs are obtained for common chaotic oscillators, such as the R\"{o}ssler system, the Lorenz system, and Chen's system.

We observe that the existence of ADIs (and thus the stabilization of fixed points) is a general characteristic of delay-coupled chaotic oscillator networks. The fact that the ADIs are all associated with a nonzero time delay indicates that time delay is a key and necessary element for the stabilization of AD states.
For fixed coupling delay, we found that stable AD as determined by the corresponding MSF either cannot occur or occurs for a continuous interval of coupling strength; however, when the coupling strength is fixed, stable AD exists for disconnected intervals of coupling delays. 
We also notice that the region within each island which corresponds to $I_{\lambda} = [-1,1]$ is smaller in area in each subsequent island as well becoming nonexistent as the islands shrink in area themselves. In Chen's system, we see that the only MSI we found does not contain a region corresponding to $I_{\lambda} = [-1,1]$ and that this implies some networks will never had a stable amplitude death state under some given parameters. 
The number, size, and location of the ADIs generally depend on the eigenvalues of the Jacobian matrix $DF$ as well as the coupling matrix. For each system considered herein, we observe that the first island (with lowest range of delay values) is always the largest, and proceeding islands shrink in area monotonically.
Within each ADI, the range of stability (visualized as the ``altitude") can vary significantly from one parameter combination to another, suggesting the importance of knowing the full MSF (as visualized as MSIs) for the determination of AD for general network topologies.

\begin{acknowledgments}
S.R.H. acknowledges support from the New York State Department of Education (Grant No.~C401608). J.S. acknowledges funding from the Simons Foundation (Grant No.~318812) and the Army Research Office (Grant No.~W911NF-12-1-0276).
\end{acknowledgments}

\bibliography{huddy_sun.bib}

\end{document}